\documentstyle[12pt,twoside]{amsart}
\begin{document}
\title{Nash problem for a toric pair and the minimal log-discrepancy}
\author{Shihoko Ishii}
\footnote{This research is partly supported by JSPS Grant-in-Aid (B), No 22340004 and (S),No19104001.}

\address{Department of Mathematics, Tokyo Institute of
Technology, Oh-Okayama, Meguro, Tokyo, Japan\\
shihoko@@math.titech.ac.jp}
\newcommand{\bC}{{\Bbb C}}
\newcommand{\bP}{{\Bbb P}}
\newcommand{\bZ}{{\Bbb Z}}
\newcommand{\bQ}{{\Bbb Q}}
\newcommand{\bR}{{\Bbb R}}
\newcommand{\bN}{{\Bbb N}}
\newcommand{\bA}{{\Bbb A}}
\newcommand{\bG}{{\Bbb G}}
\newcommand{\bV}{{\Bbb V}}
\newcommand{\R}{{\mathcal R}}
\newcommand{\ba}{{\bf a}}
\newcommand{\bp}{{\bf p}}
\newcommand{\be}{{\bf e}}
\newcommand{\bq}{{\bf q}}
\newcommand{\bl}{{\bf 1}}
\newcommand{\cale}{{\cal E}}
\newcommand{\Proj}{\operatorname{Proj}}
\newcommand{\Hom}{\operatorname{Hom}}
\newcommand{\ord}{\operatorname{ord}}
\newcommand{\orb}{\operatorname{orb}}
\newcommand{\st}{{\spec k[[t]]}}
\newcommand{\stm}{{\spec k[t]/(t^{m+1})}}
\newcommand{\sTm}{{\spec K[t]/(t^{m+1})}}
\newcommand{\sT}{{\spec K[[t]]}}
\newcommand{\tm}{{k[t]/(t^{m+1})}}
\newcommand{\D}{{\Delta}}
\newcommand{\mld}{\operatorname{mld}}
\let \cedilla =\c
\renewcommand{\c}[0]{{\mathbb C}}  
\let \crossedo =\o
\renewcommand{\o}[0]{{\mathcal O}} 
\newcommand{\z}[0]{{\mathbb Z}}
\newcommand{\n}[0]{{\mathbb N}}
\let \ringaccent=\r  
\renewcommand{\r}[0]{{\mathbb R}} 

\renewcommand{\a}[0]{{\frak{a}}} 

\newcommand{\h}[0]{{\mathbb H}}
\newcommand{\p}[0]{{\mathbb P}}
\newcommand{\f}[0]{{\mathbb F}}
\newcommand{\q}[0]{{\mathbb Q}}
\newcommand{\map}[0]{\dasharrow}
\newcommand{\qtq}[1]{\quad\mbox{#1}\quad}
\newcommand{\spec}[0]{\operatorname{Spec}}
\newcommand{\sing}[0]{\operatorname{Sing}}

\newcommand{\onto}[0]{\twoheadrightarrow}

\def\into{\DOTSB\lhook\joinrel\rightarrow}
\def\to {\longrightarrow}

\newtheorem{theorem}{Theorem}[section]
\newtheorem{lemma}[theorem]{Lemma}
\newtheorem{e-proposition}[theorem]{Proposition}
\newtheorem{corollary}[theorem]{Corollary}
\newtheorem{e-definition}[theorem]{Definition\rm}
\newtheorem{remark}{\it Remark\/}
\newtheorem{example}{\it Example\/}

\maketitle
\begin{abstract}

This paper formulates the Nash problem for a pair consisting of a toric variety and an invariant ideal and gives an affirmative answer  to the problem. 
We also prove that the minimal log-discrepacy is computed by a divisor corresponding to a Nash component, if the minimal log-discrepancy is finite.
On the other hand there exists a Nash component such that the corresponding divisor  has negative log-discrepancy,
if the minimal log-discrepancy is $-\infty$.
\end{abstract}

\section{Introduction}
\noindent
The Nash problem was posed by John F. Nash in his preprint (1968) which is published later as \cite{nash}.
The  problem is asking the bijectivity between the set of  Nash components and the set of essential divisors of a singular variety $X$.
The problem is answered positively for toric varieties and negatively in general \cite{ik}.
As the counter examples are of dimension greater than 3, the Nash problem is still 
open for surfaces and 3-folds.
The Nash problem for a surface is now steadily improving thanks to the work of M. Lejeune-Jalabert, A. Reguera (\cite{LR}, \cite{LR2}).
A Nash component is an irreducible component  of the family of arcs passing through the singular locus. 
So it does not depend on the existence of a resolution of the singularities of $X$, while an essential divisor is defined by using resolutions of the singularities of $X$.
The study of some examples gives us a feeling that we can get the information of  the singularities of $X$ from the information of the Nash components (notion without resolutions) 
even for the properties defined by using resolutions.

In this paper, we consider the Nash problem for a pair consisting of a variety and an ideal on the variety. Our principles are:
\begin{enumerate}
\item
For an object in the toric category, the Nash problem should hold.
\item
We should be able to see whether the singularities of the pair is log-canonical/log-terminal from information given by  the Nash components.
\end{enumerate}
(The first principle seems reasonable since we have some evidences 
\cite{perez}, \cite{ik}, \cite{i2}, \cite{i3}. The second principle is based on the observation for the counter example of the Nash problem \cite{ik}.)
We will show the principles are true for a toric pair consisting of a toric variety and an invariant ideal.
When we consider a pair, the primary problem is how to formulate the Nash problem for the pair.
Peter Petrov  formulated the Nash problem for a toric pair and gave an affirmative answer in \cite{p}.
But his Nash components do not satisfy (ii).
Our formulation of the Nash problem for a toric pair is different from his, but we use  his result for  our problem.
Our Nash components are constructed on a modified space of $X$ and this idea
suggests a direction for the Nash problem in the   general case.
\section{The Nash problem and minimal log-discrepancy}
\begin{e-definition}
  Let $X$ be a scheme over an algebraically closed field $k$. 
  An arc of $X$ is a $k$-morphism $\alpha:\spec K[[t]]\to X$, where $K\supset k$ is a field extension. 
  The space of arcs of $X$ is denoted by $X_\infty$ and the canonical projection 
  $X_\infty\to X$ is denoted by $\pi^X$.
  For a morphism $f:Y\to X$ of $k$-schemes, the induced morphism between the 
  arc spaces is denoted by $f_\infty: Y_\infty\to X_\infty$.
  One can find basic materials on the space of arcs in \cite{cr}.
\end{e-definition}
\vskip.3truecm

From now on we consider a pair $(X,Z)$ consisting of a variety $X$ over $k$ and a closed subscheme $Z\subset X$, or equivalently $(X, \a)$, where $\a$ is the defining ideal of $Z$. 
We always assume that $\sing X\subset |Z|$.

\begin{e-definition}
  A proper birational morphism $f:Y\to X$ with $Y$ smooth, such that $f_{Y\setminus
  f^{-1}(Z)}$ is an isomorphism on $X\setminus Z$ and $f^{-1}(Z)$ is of pure codimension 1 is called a $Z$-resolution.
  When $f$ satisfies the further conditions: $\a\o_Y$ is invertible and 
  $|f^{-1}(Z)|$ is of normal crossings, then it is called a log-resolution of $(X,Z)$.
  A divisor over $X$ is called $Z$-essential if it appears in every $Z$-resolution and 
  is called log-essential if it appears in every log-resolution. 
  \end{e-definition}
\vskip.3truecm
\begin{e-definition}
For a pair $(X, Z)$, let $f:Y\to X$ be a $Z$-resolution and $E_i (i=1,\ldots, r)$ be the irreducible exceptional divisors of $f$.
We say that $E_i$ is a $Z$-Nash divisor if the closure of 
$f_\infty\left((\pi^Y)^{-1}(E_i)\right)$ is an irreducible component of $(\pi^X)^{-1}(\sing X)$ and call this component a $Z$-Nash component. 
Note that among all divisors over $X$ there is a unique $Z$-Nash divisor up to birational equivalence for a
fixed $Z$-Nash component.
\end{e-definition}

\vskip.3truecm
\begin{theorem}[Petrov \cite{p}]
\label{pe}
Let $X$ be an affine toric variety and $Z$ an invariant closed subscheme.
Then the set of $Z$-Nash divisors and the set of $Z$-essential divisors coincide.
\end{theorem}
\vskip.3truecm
\begin{e-definition}
Let $(X,Z)$ be a pair with $X$ a normal $\bQ$-Gorenstein variety.
For a divisor $E$ over $X$, the log-discrepancy of $(X,Z)$  with respect to $E$ is 
$$a(E; X,Z):=\ord_E(K_{Y/X})-\ord_E(Z)+1,$$
where let $E$ appears on a normal variety $Y$ birational to $X$.
The minimal log-discrepacy of $(X,Z)$ is defined by
$$\mld(X,Z)=\inf\{a(E; X,Z)\mid E \ {\mbox divisor \ over\ }X\}.$$
Note that if $\dim X\geq 2$ and $\mld(X,Z)<0$, then $\mld(X,Z)=-\infty$.
A pair $(X,Z)$ is log-canonical (resp. log-terminal) if and only if 
$\mld(X,Z)\geq 0$ (resp. $\mld(X,Z)> 0$).
For a log-canonical pair $(X,Z)$, if $\mld(X,Z)=a(E;X,Z)$, then we say that 
$E$ computes the minimal log-discrepancy.
\end{e-definition}

The following shows that $Z$-Nash divisor does not necessarily compute the minimal log-discrepancy for $(X,Z)$.
The notation and terminologies on toric geometry are based on \cite{ful}.
\vskip.3truecm
\begin{example}
  Let $X$ be $\bA_\bC^3$ and $Z$ be defined by the ideal $\a=(x_1^dx_2, x_2^dx_3, x_3^dx_1)$.
  Then, $|Z|$ is the union of $x_i$-axes $(i=1,2,3)$.
  As a toric variety, $X$ is defined by a cone $\sigma:=\sum_{i=1}^3\bR_{\geq 0}\be_i$ in $N_{\bR}=\bR^3$, where $\be_1=(1,0,0),\be_2=(0,1,0),\be_3=(0,0,1)$.
 
 The $Z$-Nash divisors are $D_{\bp_i}$ $(i=1,2,3)$ which correspond to 
 $\bp_1=(0,1,1), \bp_2=(1,0,1),\bp_3=(1,1,0)$.
 When $d=2$, we can see that $\mld(X,Z)=0$, while $a(D_{\bp_i};X,Z)=1$ for 
$i=1,2,3$.
When $d=3$, we can see that $\mld(X,Z)=-\infty$, while $a(D_{\bp_i};X,Z)=1$ for 
$i=1,2,3$.
\end{example}
\vskip.3truecm
In order to produce  divisors which compute the minimal log-discrepancy, we need to modify $X$ into a more reasonable space.
We will see that for a toric pair $(X, Z)$, the normalized blow up of $X$ by the defining ideal 
$\a$ of $Z$ is an appropriate space.
\vskip.3truecm
\begin{e-definition}
  Let $(X,Z)$ be a toric pair and
 let $\varphi: \overline X\to X$ be the normalized blow up by the defining ideal $\a$ 
 of $Z$.
  Let $f:Y\to X$ be a log-resolution and $E_i (i=1,\ldots, r)$ be the irreducible exceptional divisors of $f$, then $f$ factors as $f=\varphi\circ g$ for $g:Y\to \overline X$.
We say that $E_i$ is a log-Nash divisor for  $(X, Z)$ if the closure of 
$g_\infty\left((\pi^Y)^{-1}(E_i)\right)$ is an irreducible component of $(\pi^{\overline X})^{-1}(\varphi^{-1}(Z))$ and call this component a log-Nash component. 
Note that among all divisors over $X$ there is a unique log-Nash divisor up to birational equivalence for a fixed log-Nash component.
\end{e-definition}
\vskip.3truecm
\begin{theorem}
  Let $(X,Z)$ be a toric pair, then the following hold:
\begin{enumerate}
\item The set of log-Nash divisors for $(X,Z)$ coincides with the set of log-essential divisors.
\item If $X$ is $\bQ$-Gorenstein and $(X,Z)$ is log-canonical, then a  log-Nash divisor computes the minimal log-discrepancy.
\item If $X$ is $\bQ$-Gorenstein and $(X,Z)$ is not log-canonical, then there is a 
log-Nash divisor  with a negative log-discrepancy.
\end{enumerate}
\end{theorem}

\begin{pf}
 First of all, note that the normalized blow up $\varphi: \overline X\to X$ is  a toric morphism. 
  Actually it corresponds to the decomposition into the dual fan of the
  Newton polygon $\Gamma_+(\a)$ for the ideal $\a$ of $Z$.
  In the same way as in \cite[Theorem 2.15]{ik}, we have the first inclusion of the following, while the second one is trivial:

 \ \ \ \ \ \{log-Nash divisors for $(X,Z)$\}$\subset$ \{log-essential divisors for $(X,Z)$\}.
  
 \ \ \ \ \ \ \ \ \ \ \ \ \ \ \ \ \ \ \ \ \ \ \ \ \ \ \ \ \ \ \ \ \ \ \ \ \ \ \ \ \ \ \ \  $\subset$ \{divisors appearing in every toric log-resolution of $(X,Z)$\}. 
\newline 
For the statement (i), it is sufficient to show the equality of the first and the third sets.
In fact, the first set is the same as \{$\varphi^{-1}(Z)$-Nash divisors\}   and
it coincides with \{$\varphi^{-1}(Z)$-essential divisors\} by Petrov's result
Theorem\ref{pe}.
His proof also shows that this set is the same as \{divisors appearing in every toric
$\varphi^{-1}(Z)$-resolution\}.
  As a toric log-resolution always factor through $\overline X$ and an invariant divisor on a non-singular toric variety is always normal crossings, 
  therefore a toric 
  $\varphi^{-1}(Z)$-resolution of  $\overline X$ is the same as a toric 
  log-resolution of $(X,Z)$.
  Thus, it follows the required coincidence of the sets.

  In order to prove (ii) and (iii), we remark that for a prime divisor $E$ and an effective divisor $D$ on a 
  non-singular variety $Y$ and the generic point $\gamma$ of  $ (\pi^Y)^{-1}(E)$,
  $$\ord_E(D)=\ord_\gamma(\o_Y(-D)).$$
  Let $r$ be the index of $K_X$.
  Let ${\mathcal L}=\varphi^*\omega_X^{[r]}\otimes_{\o_{\overline X}}{\o_{\overline X}}
  (r\varphi^{-1}(Z))$, then it is an invertible sheaf on ${\o_{\overline X}}$.
  If $(X,Z)$ is log-canonical, then ${\mathcal L}$ is moreover an ideal sheaf of ${\o_{\overline X}}$.
   Indeed, on a toric log-resolution $f:Y\to X$, we have

$$a(E;X,Z)=\frac{1}{r}\left(\ord _E (rK_{Y/X})-r\ord_E f^{-1}(Z) +r\right)$$
\begin{equation}
\label{a}
=\frac{1}{r}\ord_E(-f^*(rK_X)-rf^{-1}(Z))
\geq 0
\end{equation}
   for every ivariant prime divisor $E$ on $Y$.
   Here we used that $$K_Y=-\sum _{D:invariant\ prime\ divisor}D.$$
   Therefore, $-f^*(rK_X)-rf^{-1}(Z)$ is an effective Cartier divisor on $Y$, which yields
   that ${\mathcal L}=g_*g^*({\mathcal L})=g_*\left(\o_Y(f^*(rK_X)+rf^{-1}(Z))\right)\subset 
   \o_{\overline X}$, where $g:Y\to \overline X$ is the factorization of $f$ by 
   $\varphi:\overline X\to X$.
   Now, we see by (\ref{a}), $$a(E;X,Z)=
   \frac{1}{r}\ord_\gamma(g^*{\mathcal L})=\frac{1}{r}\ord_{g_\infty(\gamma)}({\mathcal L}),$$
   where $\gamma$ is the generic point of  $ (\pi^Y)^{-1}(E)$. 
   If $E$ computes the minimal log-discrepancy, take the log-Nash divisor $E_0$ with the generic point $\alpha\in (\pi^Y)^{-1}(E_0)$ such that 
   $g_\infty(\gamma)\in \overline{g_\infty(\alpha)}$.
   Then, by the upper semi-continuity of the order in the arc space,
   $$\mld(X,Z)=\ord_{g_\infty(\gamma)}({\mathcal L})\geq \ord_{g_\infty(\alpha)}({\mathcal L})=a(E_0;X,Z),$$
   which shows that $E_0$ computes the minimal log-discrepancy as required in (ii).
  
  If $(X,Z)$ is not log-canonical, then ${\mathcal L}\not\subset \o_{\overline X}$.
  Indeed, if   ${\mathcal L}\subset \o_{\overline X}$, for every prime divisor $E$ over $X$
   with the generic point $\gamma\in (\pi^Y)^{-1}(E)$, we have 
   $a(E;X,Z)=\frac{1}{r}\ord_\gamma \o_Y(g^*{\mathcal L})\geq 0$, which implies
   that $(X,Z)$ is log-canonical, a contradiction.
   Now, we can put  ${\mathcal L}=\o_{\overline X}(D-D')$, where $D>0$ and $D'\geq 0$ are invariant divisors which do not have common components.
   Then, for a prime divisor $E\leq D$
   $$a(E;X,Z)=\frac{1}{r}\ord_E(-D+D') <0.$$
   As $E$ is a prime divisor on $\overline X$ with the support in $|\varphi^{-1}(Z)|$, it is a log-Nash divisor.
\end{pf}

\begin{remark}
  One can prove (ii) and (iii) also by the combinatorics on the fan.
  
  Another way to prove (ii) is  to observe that every divisor that computes
  the minimal log-discrepancy is an log-essential divisor and then 
  use (i).
  This way  provides us with the stronger  fact that all divisors that compute the minimal log
  discrepancy are  log-Nash divisors.
  But I presented a proof above which does not use (i), because this proof may be useful to study  a general case in which the Nash problem does not hold.
\end{remark}






\end{document}